\documentclass[11pt,a4paper]{article}
\usepackage{amsfonts,amscd,amsthm,amsgen,amsmath,amssymb}
\usepackage[all]{xy}
\newtheorem{prop}{Proposition}[section]
\newtheorem{thm}{Theorem}[section]

\newtheorem{lemp}{Lemma}[section]

\newtheorem{corp}{Corollary}[section]
\theoremstyle{remark}
\newtheorem{rem}{Remark}[section]

\theoremstyle{definition}
\newtheorem{defn}{Definition}[section]
\newcommand{\leftexp}[2]{{\vphantom{#2}}^{#1}{#2}}
\title{Moduli of Coassociative Submanifolds and Semi-Flat Coassociative Fibrations}
\usepackage{amsfonts,amscd,amsthm,amsgen,amsmath,amssymb}
\author{David Baraglia}
\date{\today}
\begin{document}
\maketitle
\begin{abstract}
We show that the moduli space of deformations of a compact coassociative submanifold $C$ has a natural local embedding as a submanifold of $H^2(C,\mathbb{R})$. We show that a $G_2$-manifold with a $T^4$-action of isomorphisms such that the orbits are coassociative tori is locally equivalent to a minimal $3$-manifold in $\mathbb{R}^{3,3}$ with positive induced metric where $R^{3,3}\cong H^2(T^4,\mathbb{R})$. By studying minimal surfaces in quadrics we show how to construct minimal $3$-manifold cones in $\mathbb{R}^{3,3}$ and hence $G_2$-metrics from equations similar to a set of affine Toda equations.
The relation to semi-flat special Lagrangian fibrations and the Monge-Amp\`ere equation are explained.
\end{abstract}

\section{Introduction}
The well-known conjecture of Strominger, Yau and Zaslow \cite{syz} provides a geometric picture of mirror symmetry, at least in the so-called large complex structure limit. The conjecture proposes mirror pairs of Calabi-Yau manifolds which are special Lagrangian torus fibrations over the same base, but with dual fibres. In understanding limiting cases of the conjecture one is motivated to study {\em semi-flat} special Lagrangian fibrations. These are fibrations in which the fibres are flat tori. It is known in this case that the base has natural affine coordinates and a function $\phi$ satisfying the real Monge-Amp\`ere equation ${\rm det}H(\phi)= 1$, where $H(\phi)$ is the Hessian of $\phi$ with respect to the affine coordinates \cite{hit}.\\

In M-theory $G_2$-manifolds play a role equivalent to Calabi-Yau manifolds in string theory, so it is natural to ask whether there is an analogue of the SYZ conjecture for $G_2$-manifolds. Gukov, Yau and Zaslow \cite{gyz} argue that $G_2$-manifolds admitting Calabi-Yau string theory duals are fibred by coassociative submanifolds. Therefore it seems natural to expect a corresponding mirror symmetry relating pairs of coassociative fibrations over the same base. We take this as motivation to study $G_2$-manifolds fibred by flat coassociative tori. More specifically we call a coassociative fibration $X$ {\em semi-flat} if there is a $T^4$-action of isomorphisms of $X$ such that the orbits are coassociative submanifolds. The key result is Theorem \ref{thethm} which states that the base $B$ locally maps into $H^2(T^4,\mathbb{R})$ (equipped with the intersection form) as a positive definite minimal $3$-submanifold and conversely such a minimal submanifold gives rise to a semi-flat coassociative fibration.\\

We investigate the structure of the moduli space of deformations of a compact coassociative submanifold. Adapting the approach of \cite{hit} which studies the moduli space of special Lagrangians, we find that the moduli space $\mathcal{M}$ of deformations of a compact coassociative submanifold $C$ has locally a natural map $u : \mathcal{M} \to H^2(C,\mathbb{R})$ defined up to an affine map. The $L^2$ metric on the moduli space is then the induced metric under $u$.\\

Finding examples of compact coassociative fibrations other than from special Lagrangian fibrations is a very difficult task. We are however able to show that for a $G_2$-manifold with finite fundamental group, a coassociative fibration must have singularities. To understand the types of singularities that one might expect to encounter we show how to construct examples of compact coassociative fibrations for a $G_2$-structure with torsion. More specifically we construct $7$-manifolds with $G_2$-structures such that the corresponding $4$-form is closed but the $3$-form need not be. These manifolds have coassociative fibrations degenerating over a smooth link in $S^3$. Moreover the possible singular fibres (except for a finite number of fibres) can be readily determined from \cite{mat}.\\

We consider the case of a $G_2$-manifold that is a product $X = Y \times S^1$ of a Calabi-Yau manifold $Y$ and a circle. In this case the coassociative and special Lagrangian moduli spaces are related. We show that in the semi-flat case the minimal submanifold equations reduce to the Monge-Amp\`ere equation, recovering the known result on semi-flat Calabi-Yau manifolds.\\

In \cite{lyz} the authors seek solutions to the Monge-Amp\`ere equation on a $3$-dimensional base that is a cone. This reduces to equations on a surface, in fact the equations for an elliptic affine sphere. This amounts to solving the following equation
\begin{equation*}
\psi_{z\overline{z}} + |U|^2e^{-2\psi} + \tfrac{1}{2}e^\psi = 0
\end{equation*}
where $U$ is a holomorphic cubic differential. This is a real form of the affine Toda field equations for the affine Dynkin diagram $A^{(2)}_2$, studied by Tzitz\'eica \cite{tzitz}. In a similar fashion we seek to reduce from the minimal submanifold equations on a $3$-manifold to equations on a surface. We consider semi-flat $G_2$-manifolds with a vector field commuting with the $T^4$-action which essentially scales the associative $3$-form. This corresponds to the minimal $3$-manifold being a cone, which in turn is equivalent to a minimal surface in the quadric of unit vectors.\\

We outline the general situation of minimal surfaces in quadrics, adapting to indefinite signature the theory of harmonic sequences. We show that a class of minimal surfaces into odd dimensional quadrics called superconformal are equivalent to a set of equations that are a real form of the $D^{(2)}_{r+1}$ affine Toda equations. In the case of the unit quadric in $\mathbb{R}^{3,3}$ the equations obtained are
\begin{eqnarray*}
2(w_1)_{z\overline{z}} &=& -e^{2w_2-2w_1} - e^{2w_1}, \\
2(w_2)_{z\overline{z}} &=& q\overline{q}e^{-2w_2} + e^{2w_2 - 2w_1}.
\end{eqnarray*}
where $q$ is a holomorphic cubic differential. The case where the $G_2$-manifold is a product of a Calabi-Yau manifold and a circle corresponds to the reduction $e^{2w_2} = q\overline{q}e^{-2w_1}$ in which case the equations reduce to the equation of Tzitz\'eica.\\


\section{Deformations of coassociative submanifolds}\label{seccoasub}
Let $X$ be a $G_2$-manifold, that is a Riemannian manifold with holonomy a subgroup of $G_2$. We will be concerned with the moduli space of deformations of a coassociative submanifold $C$ of $X$.\\

Recall that $G_2$ is the stabiliser of a $3$-form $\phi_0$ on $\mathbb{R}^7$ and is a subgroup of ${\rm SO}(7)$ so also stabilises the $4$-form $\psi_0 = *\phi_0$. We use the convention in \cite{joyce} for the standard $3$ and $4$-forms on $\mathbb{R}^7$. Hence if ${e^1, \dots, e^7}$ are a standard basis of $1$-forms for $\mathbb{R}^7$, the standard $3$ and $4$-forms $\phi_0,\psi_0$ are taken as
\begin{equation}
\begin{aligned}
\phi_0 &= e^{123} + e^1 \! \! \wedge \! (e^{45} + e^{67}) + e^2 \! \! \wedge \! (e^{46} - e^{57}) + e^3 \! \! \wedge \! (-e^{47} - e^{56}) \\
\psi_0 &= e^{4567} + e^{23} \! \! \wedge \! (e^{45} + e^{67}) + e^{31} \! \! \wedge \! (e^{46} - e^{57}) + e^{12} \! \! \wedge \! (-e^{47} - e^{56})
\end{aligned}
\end{equation}
where we use the notation $e^{ij\dots k}=e^i \wedge e^j \wedge \dots \wedge e^k$.\\

A $G_2$-manifold $X$ therefore is a Riemannian $7$-manifold with a reduction of structure from $O(7)$ to $G_2$ such that the corresponding $3$-form $\phi$ is covariantly constant. An equivalent characterisation is that $X$ is a $7$-manifold with a reduction of structure to $G_2$ such that the $3$-form $\phi$ and $4$-form $\psi = *\phi$ are closed \cite{joyce}.\\

Recall that a $4$-dimensional submanifold $C$ in $X$ is {\em coassociative} if and only if $\phi|_C=0$ and in such a case $C$ inherits an orientation $\psi|_C$. We can relax the submanifold condition to just the requirement that $f : C \to X$ is an immersion in which case the restriction notation $\phi|_C$ should be understood as the pullback $f^*\phi$. The normal bundle $NC$ of a coassociative submanifold $C$ can be identified with $\wedge^2_+ T^*C$, the bundle of self-dual $2$-forms on $C$ by sending a normal vector $\nu$ to the $2$-form $\iota_\nu \phi |_C$.\\

Two immersions $f,g : C \to X$ are called {\em isotopic} if there is a map $F : C \times [0,1] \to X$ such that each map $f_t = F( \; , t) : C \to X$ is an immersion and $f = f_0, g = f_1$. Moreover if $f$ and $g$ are coassociative submanifolds then we say $F$ is a {\em deformation} through coassociative submanifolds if $f_t^*\phi = 0$ for all $t$. Note that if $C$ is compact then by means of a $t$-dependent reparametrisation we may assume that $\partial_t F$ is a normal vector field for all $t$.\\

Given a compact coassociative submanifold $C$ let $\mathcal{M}$ denote the moduli space of all coassociative submanifolds that are deformations of $C$. McLean \cite{mclean} proves the following:
\begin{prop}\label{mclean1}
A normal vector field $\nu$ to a compact coassociative submanifold $C$ is the vector field normal to a deformation through coassociative submanifolds if and only if the corresponding $2$-form $\iota_\nu \phi|_C$ is a closed self-dual $2$-form, hence harmonic. There are no obstructions to extending a first order deformation to a family of actual deformations.
\end{prop}

Actually, McLean uses a convention for the associative $3$-form leading to anti-self dual harmonic forms. It follows from the work of McLean that $\mathcal{M}$ is a smooth manifold of dimension $b^2_+(C)$. Moreover there is a natural fibre bundle $\pi : \mathcal{C} \to \mathcal{M}$ with fibres diffeomorphic to $C$ and map $e: \mathcal{C} \to X$ such that for each $m \in \mathcal{M}$, $e|_{\pi^{-1}(m)} : \pi^{-1}(m) \to X$ is the coassociative submanifold corresponding to $m \in \mathcal{M}$. We let $C_m$ denote $\pi^{-1}(m)$ and $e_m : C_m \to X$ the corresponding immersion.\\

From the proposition have for each $m\in \mathcal{M}$ a natural isomorphism $\lambda_m : T_m\mathcal{M} \to \mathcal{H}^2_+(C_m,\mathbb{R})$, where $\mathcal{H}^2_+(C_m,\mathbb{R})$ is the space of harmonic self-dual $2$-forms on $C_m$. Since the self-dual cohomology on $C_m$ depends on the metric which varies with $m$ we do not have a natural identification of all the $\mathcal{H}^2_+(C_m,\mathbb{R})$. However we can still take the cohomology class $[\lambda_m] \in H^2(C_m , \mathbb{R})$. Working locally, we will not distinguish between $\mathcal{M}$ and a neighborhood of $\mathcal{M}$ over which the fibre bundle $\mathcal{C}$ is trivial. It follows that we can canonically identify the cohomology of the fibres. Fixing a single fibre $C = C_{m_0}$ we can identify cohomology of each fibre with $H^2(C,\mathbb{R})$. Under these identifications $\lambda_m$ determines a map $[\lambda] : T\mathcal{M} \to H^2(C,\mathbb{R})$, that is a $H^2(C,\mathbb{R})$-valued $1$-form. We find

\begin{prop}\label{closed}
The $1$-form $[\lambda]$ on $\mathcal{M}$ is closed. Thus locally we may write $[\lambda] = du$ where $u : \mathcal{M} \to H^2(C,\mathbb{R})$ is an immersion unique up to a translation.
\begin{proof}
Recall that we have a fibre bundle $\pi : \mathcal{C} \to \mathcal{M}$ and a natural map $e : \mathcal{C} \to X$ representing the full local family of deformations. Fix a basepoint $m_0 \in \mathcal{M}$ and let $C = C_{m_0}$. Since we are working locally in $\mathcal{M}$ we assume a local trivialisation $\mathcal{C} = \mathcal{M} \times C$.\\

Since $C$ is a calibrated submanifold it is oriented and we have also assumed $C$ is compact, hence given a homology class $A \in H_2(C,\mathbb{R})$ let $\eta_A \in \Omega^{2}(C,\mathbb{R})$ represent the Poincar\'e dual, that is for any $\sigma \in H^2(C,\mathbb{R})$,
\begin{equation*}
\langle \sigma , A \rangle = \int_C \sigma \wedge \eta_A.
\end{equation*}
Extend $\eta_A$ to a form on $\mathcal{C} = \mathcal{M} \times C$ in the natural way. We claim that the $2$-form $\langle \lambda , A \rangle$ is given by $p_* (e^*\phi \wedge \eta_A)$ where $p_*$ represents integration over the fibres of $\mathcal{C}$. Since $\phi$ is closed and integration over fibres takes closed forms to closed forms, this will prove the proposition.\\

Let $t^1,\dots ,t^m$ be local coordinates on $\mathcal{M}$. Then by the definition of $e$ we have for each $i$ and each $m \in \mathcal{M}$ a section $\nu_i$ of the bundle $e_m^*(TX)$ representing the corresponding deformation, given by $\nu_i(c) = e_*(m,c)\partial /\partial t^i$. Note that $\nu_i$ need not be normal to $C_m$. Define the $2$-forms $\theta_i$ on $C_m$ by $\theta_i = (\iota_{\nu_i} \phi)|_{C_m}$. Since each submanifold $C_m$ is coassociative we have $\phi|_{C_m}=0$. Thus only the normal component of $\nu_i$ contributes to $\theta_i$ so the $\theta_i$ are harmonic by Proposition \ref{mclean1}. Moreover since $\mathcal{M}$ represents all local deformations the $\theta_i$ span $\mathcal{H}^2_+(C_m,\mathbb{R})$. We also have $[\lambda] = [\theta_i] dt^i$.\\

Now consider $e^*\phi$. Since $(e^*\phi)|_{C_m}=0$ it follows that we can write $e^*\phi = dt^i \wedge \hat{\theta}_i$ for some $2$-forms $\hat{\theta}_i$. Moreover we see that $\hat{\theta}_i|_{C_m} = \theta_i(m)$. Now when we perform integration over the fibres of $\mathcal{C}$ we find
\begin{eqnarray*}
p_* (e^*\phi \wedge \eta_A) &=& p_* ( dt^i \wedge \hat{\theta}_i \wedge \eta_A) \\
&=& dt^i \int_{C_m} \hat{\theta_i}|_{C_m} \wedge \eta_A \\
&=& dt^i \int_A \hat{\theta_i}|_{C_m} \\
&=& dt^i\int_A \theta_i \\
&=& \langle \lambda , A \rangle.
\end{eqnarray*}
\end{proof}
\end{prop}

The moduli space $\mathcal{M}$ has a natural metric $g_{L^2}$ which we call the $L^2$ metric. For $X,Y \in T_m\mathcal{M}$ the metric is defined as
\begin{equation}
g_{L^2}(X,Y) = \int_{C_m} \lambda_m(X) \wedge \lambda_m(Y).
\end{equation}
From Proposition \ref{closed} we immediately have
\begin{prop}
The natural $L^2$ metric $g_{L^2}$ on $\mathcal{M}$ is the pull-back under $u:\mathcal{M} \to H^2(C,\mathbb{R})$ of the non-degenerate inner product on $H^2(C,\mathbb{R})$ given by the intersection form. Since $u$ maps each tangent space to a positive definite subspace, $g_{L^2}$ is positive definite.
\end{prop}

Next we turn to the issue of monodromy. Suppose we have two points $s,t \in \mathcal{M}$ and two curves joining them. Along the curves we get isotopies between $C_s$ and $C_t$ and hence corresponding isomorphisms of their cohomology. If the two curves are not homotopic these two isomorphisms need not agree, though in either case a basis for $H^2(C_s,\mathbb{Z})$ (modulo torsion) will be sent to a basis for $H^2(C_t,\mathbb{Z})$ preserving the intersection form. Therefore the ambiguity is an element of ${\rm \bf SO}(Q,\mathbb{Z}) = {\rm \bf SO}(Q) \cap {\rm \bf SL}(b^2,\mathbb{Z})$ where $Q$ denotes the intersection form and $b^2$ is the second Betti number of $C$. The unit determinant condition follows since for any closed curve we have an isotopy through coassociative submanifolds preserving the orientation induced by the calibrating form $\psi$. In the case where $C = T^4$, cohomology is generated by the $1$-cocycles so the ambiguity can be thought of as an element of ${\rm \bf SL}(4,\mathbb{Z})$ acting via the representation ${\rm \bf SL}(4,\mathbb{Z}) \to {\rm \bf SO}(3,3,\mathbb{Z})$ on $2$-forms.\\

The above ambiguity prevents us from defining $[\lambda]$ as a global $1$-form on $\mathcal{M}$. However we can define $[\lambda]$ on the universal cover $\hat{\mathcal{M}}$ of $\mathcal{M}$. Then since $[\lambda]$ is closed and $\hat{\mathcal{M}}$ simply connected we may write $[\lambda] = d\hat{u}$ giving a developing map $\hat{u} : \hat{\mathcal{M}} \to H^2(C,\mathbb{R})$ of the universal cover. Note that while $[\lambda]$ is defined on $\mathcal{M}$ up to ambiguity in ${\rm \bf SO}(Q,\mathbb{Z})$, there is a second ambiguity in $u$, namely translations. Therefore the monodromy representation has the form $\rho : \pi_1(\mathcal{M}) \to {\rm \bf SO}(Q,\mathbb{Z}) \ltimes \mathbb{R}^{b^2(C)}$.

\section{Global aspects of coassociative fibrations}\label{coassf}
We now turn our attention to the study of coassociative fibrations. More specifically we look at some global aspects of coassociative fibrations of compact manifolds.\\

We will be considering fibrations $\pi : X \to B$ of $G_2$-manifolds such that the fibres are coassociative submanifolds. We can also consider fibrations which degenerate, that is $\pi$ need not be a locally trivial fibre bundle. In fact we will show that if $X$ is compact and has finite fundamental group then any fibration coassociative or not must degenerate.

\begin{prop}
Let $X$ be a compact $G_2$-manifold with finite fundamental group. Then $X$ admits no locally trivial fibre bundles $X \to B$ onto a $3$-dimensional base.
\begin{rem}
A compact $G_2$-manifold $X$ has finite fundamental group if and only if the holonomy of $X$ is the whole of $G_2$ \cite{joyce}.
\end{rem}
\begin{proof}
Our argument expands upon \cite{cavalcanti} ex 8.4, pg 130. We may assume $X$ is connected. Further it suffices to replace $X$ by its universal cover which is also compact, so we assume $X$ is simply connected. Let $\pi : X \to B$ be a locally trivial fibre bundle where $B$ is a $3$-manifold and let $F$ be one of the fibres. If $F$ is not connected we may replace $B$ by its universal cover $\tilde{B}$ to get another fibration $\pi_1 : X \to \tilde{B}$ with fibre equal to a connected component of $F$. So we may as well assume $B$ is simply connected and $F$ is connected. 

Now since $B$ is a compact simply connected $3$-manifold it is oriented and $H_1(B,\mathbb{Z}) = H_2(B,\mathbb{Z}) = 0$. Then by the Hurewicz theorem $\pi_2(B) = 0$ also (of course the Poincar\'e conjecture implies the base is diffeomorphic to the $3$-sphere but we don't need this fact). The long exact sequence of homotopy groups implies that $\pi_1(F) = 0$. Further $F$ must be oriented since $X$ is.\\

We will make use of the Leray-Serre spectral sequence (with coefficients in $\mathbb{R}$) in order to gain information on the cohomology of $X$. Since $B$ is simply connected we have $E^{p,q}_2 = H^p(B,\mathbb{R}) \otimes H^q(F,\mathbb{R})$. In order to calculate $H^2(X,\mathbb{R})$ and $H^3(X,\mathbb{R})$ there is only one relevant non-trivial differential to consider $d_3 : E^{0,2}_3 \to E^{3,0}_3$, with $E^{0,2}_3 = H^2(F,\mathbb{R})$ and $E^{3,0}_3 = H^3(B,\mathbb{R})$. We then have $H^2(X,\mathbb{R}) = {\rm ker}(d_3)$ and $H^3(X,\mathbb{R}) = {\rm coker}(d_3)$. However the differential $d_3$ must vanish for otherwise we have $H^3(X,\mathbb{R}) = 0$ which is impossible on a compact $G_2$-manifold. So $d_3 = 0$ and it follows that the following maps are isomorphisms:
\begin{equation}
\begin{aligned}
i^* &: H^2(X , \mathbb{R}) \to H^2(F,\mathbb{R}) \\
\pi^* &: H^3(B,\mathbb{R}) \to H^3(X,\mathbb{R})
\end{aligned}
\end{equation}
where $i^*$ is induced from the inclusion of some fibre $i : F \to X$ and $\pi^*$ is induced by the projection $\pi : X \to B$.\\

Let $\phi$ be the $G_2$ $3$-form on $X$. Then the cohomology class of $\phi$ has the form $\left[ \phi \right] = c \pi^* \left[ dvol_B \right]$ where $dvol_B$ is a volume form on $B$ such that $\int_B dvol_B = 1$ and $c$ is some non-vanishing constant. We claim that $\pi^*(dvol_B)$ is Poincar\'e dual to the fibre $F$. This is a straightforward consequence of fibre integration. Therefore if $\mu$ is a closed $4$-form on $X$ we have
\begin{equation}\label{poinc}
\int_F i^*\mu = \int_X \mu \wedge \pi^*(dvol_B) = c^{-1} \int_X \mu \wedge \phi.
\end{equation}
Let us recall two cohomological properties of compact $G_2$-manifolds \cite{joyce} where we continue to assume that $H^1(X,\mathbb{R}) = 0$. First there is a symmetric bilinear form $\langle \, , \, \rangle$ on $H^2(X,\mathbb{R})$ given by
\begin{equation}\label{pairing}
\langle \eta , \xi \rangle = \int_X \eta \wedge \xi \wedge \phi.
\end{equation}
This form is negative definite. Secondly if $p_1(X) \in H^4(X,\mathbb{R})$ is the first Pontryagin class of $X$ then
\begin{equation}\label{lessthanzero}
\int_X p_1(X) \wedge \phi \, < 0.
\end{equation}
Combining (\ref{pairing}) with (\ref{poinc}) and the fact that $i^* : H^2(X,\mathbb{R}) \to H^2(F,\mathbb{R})$ is an isomorphism we find that the intersection form on $F$ is negative definite. Hence Donaldson's theorem (\cite{don}, Theorem 1.3.1) implies the intersection form of $F$ is diagonalisable, i.e. of the form ${\rm diag}(-1,-1, \dots , -1)$. Now $F$ is a spin manifold because $X$ is spin and the normal bundle $N$ of $F$ in $X$ is trivial so that $w_2(TF) = w_2(TF \oplus N) = i^* w_2(TX) = 0$ where $w_2$ denotes the second Stiefel-Whitney class. But now from Wu's formula the intersection form of $F$ must be even. Therefore the intersection form must be trivial and $H^2(F,\mathbb{R}) = 0$.

Now we also have that $p_1(F) = p_1(TF \oplus N) = i^*p_1(X)$. Therefore
\begin{equation*}
\int_X p_1(X) \wedge \phi = c \int_F p_1(F) = 0
\end{equation*}
where the last equality follows from the Hirzebruch signature theorem. But this contradicts (\ref{lessthanzero}), hence such fibre bundles $\pi : X \to B$ can not exist.
\end{proof}
\end{prop}

Next we consider the question of what smooth fibres a compact coassociative fibration can have. Clearly if $F$ is such a fibre we must have $b^2_+(F) \ge 3$ and moreover $\wedge^2_+ T^*F$ has a trivialisation by harmonic forms (in particular $F$ has an ${\rm SU}(2)$-structure). We also have:
\begin{prop}\label{riemfibr}
Let $\pi : X \to B$ be a coassociative fibration with compact fibres. Then the base $B$ can be given a metric such that $\pi$ is a Riemannian submersion (around non-singular fibres) if and only if the fibres are Hyperk\"ahler. Moreover in this case the base identifies with the moduli space of deformations of a fibre and the base metric $g_B$ is related to the moduli space $L^2$ metric $g_{L^2}$ by
\begin{equation}\label{basemetric}
g_B = \frac{1}{2 {\rm vol(F)}}g_{L^2}.
\end{equation}
\begin{proof}
Let $b \in B$ and $F = \pi^{-1}(b)$. Choose a basis $v_i$ for $T_bB$ and let $\tilde{v}_i$ be the horizontal lifts. Define corresponding harmonic $2$-forms $\omega_i$ by $\omega_i = \iota_{\tilde{v}_i} \phi |_F$. Then the $\omega_i$ are a frame for $\wedge^2_+T^*F$. Now if $g$ is the metric on $X$ then
\begin{equation}\label{gdvol1}
g(A , B) dvol_X = \frac{1}{6}\phi \wedge \iota_A \phi \wedge \iota_B \phi.
\end{equation}
From this it follows that
\begin{equation}\label{wedge}
\omega_i \wedge \omega_j = 2g(\tilde{v}_i , \tilde{v}_j ) dvol_F.
\end{equation}
Now suppose $B$ has a metric such that $\pi$ is a Riemannian submersion. Choose the $v_i$ to be an orthonormal basis. Then the $\omega_i$ are Hyperk\"ahler forms on $F$. Conversely suppose $F$ is a compact Hyperk\"ahler $4$-manifold. Then since $b^2_+(F) = 3$ we have that the space of self-dual harmonic $2$-forms is $3$-dimensional, spanned by the Hyperk\"ahler forms. Thus the $\omega_i$ are constant linear combinations of the Hyperk\"ahler forms, hence $g(\tilde{v}_i , \tilde{v}_j)$ is constant along the fibres. So $B$ can be given a metric $g_B$ making $\pi$ a Riemannian submersion. In this case since $b^2_+(F) = 3$ we see that the base exhausts all deformations of a fibre through coassociative submanifolds so that the base identifies with the moduli space of deformations. Moreover if we integrate (\ref{wedge}) over $F$ we get (\ref{basemetric}).
\end{proof}
\end{prop}
Note that a smooth Hyperk\"ahler $4$-manifold is either a torus or a K3 surface. Therefore these are likely candidates for the fibres of a compact coassociative fibration.\\

We now move on to the question of what sort of singularities can occur for a compact coassociative fibration $f:X \to B$ and what does the discriminant locus $\Delta = \{ b \in B \, | \, \exists \, x \in f^{-1}(b), \, {\rm rank}(df_x) < 3 \} \subset B$ look like? 

In the case of special Lagrangian fibrations $Y \to B$ that are sufficiently well behaved the discriminant locus has codimension $2$ \cite{bai}. Under assumptions that the singularities are well behaved Baier \cite{bai} shows that if $\Delta$ is smooth then $Y$ has vanishing Euler characteristic, so generally we expect $\Delta$ not to be smooth.

Returning to the case of coassociative fibrations $X \to B$ we might likewise expect under reasonable assumptions on the singularities that the discriminant locus $\Delta \subset B$ has codimension $2$. Interestingly the constraint on smoothness in the special Lagrangian case no longer seems to be an issue since $X$ is odd dimensional, hence always has vanishing Euler characteristic. Moreover Kovalev has constructed examples of compact coassociative K3 fibrations with discriminant locus a smooth link \cite{kov}.\\

We will not investigate the issue of singularities in any depth and instead we will simply provide a model for producing examples of compact coassociative fibrations on manifolds with $G_2$-structures with torsion, that is the $3$-form will not be closed.

Suppose we have a Hyperk\"ahler $8$-manifold $M$ with holomorphic Lagrangian fibration $\pi : M \to \mathbb{CP}^2$ (with singularities). That is if $\omega_I, \omega_J, \omega_K$ are the Hyperk\"ahler forms then the non-singular fibres of $\pi$ are complex submanifolds with respect to $\omega_I$ and are Lagrangian with respect to the holomorphic symplectic form $\Omega = \omega_J + i\omega_K$. Note that the non-singular fibres are necessarily tori. We can give $M$ the structure of a ${\rm Spin}(7)$-manifold where the $4$-form is
\begin{equation}
\Phi = \frac{1}{2}\omega_I^2 + \frac{1}{2}\omega_J^2 - \frac{1}{2}\omega_K^2.
\end{equation}
This makes the fibres of $\pi$ into Cayley $4$-folds, indeed since we can write $\Phi = \omega_I^2/2 + \Omega \wedge \overline{\Omega}/2$ it follows that $\Phi|_{{\rm ker}\pi_*} = \omega_I^2/2$ which is the volume form on the fibres since the fibres are complex submanifolds. We have therefore produced an example of a compact Cayley $4$-fold fibration (albeit with holonomy in ${\rm Sp}(2)$). Note that the discriminant locus $\Delta \subset \mathbb{CP}^2$ is an algebraic curve in $\mathbb{CP}^2$.

To get a coassociative fibration let us take a smooth embedded $3$-sphere $S^3 \subset \mathbb{CP}^2$ that meets $\Delta$ transversally and avoids any singular points of $\Delta$. Therefore if $\Delta$ has real dimension $2$ then $\Delta \cap S^3$ will be a smooth link in $S^3$. Let $X = \pi^{-1}(S^3)$. We claim that $X$ has an almost $G_2$-structure such that $\pi: X \to S^3$ is a coassociative fibration with discriminant locus the link $\Delta \cap S^3$. To see this let $v$ be the unit normal to $TX$ and let $v^* = g(v, \, )$ where $g$ is the metric on $M$. Then on $TM|_X$ we may write
\begin{equation}
\Phi = v^* \wedge \phi + \psi
\end{equation}
where $\iota_v \phi =  \iota_v \psi = 0$. Then $\phi$ defines an associative $3$-form on $X$ which is generally not closed and $\psi$ is the corresponding $4$-form. Note on the other hand that $\psi$ is closed since $\Phi|_X = \psi$. This also shows that the smooth fibres of $\pi$ are coassociative submanifolds.
To make this example more explicit let $f:S \to \mathbb{CP}^1$ be an elliptic fibration of a K3 surface over $\mathbb{CP}^1$. Then we will take our Hyperk\"ahler $8$-manifold $M$ to be the Hilbert scheme ${\rm Hilb}^2S$ of pairs of points on $S$. This can be concretely described as follows: the product $S \times S$ has a natural $\mathbb{Z}_2$-action given by exchanging points. The Hilbert scheme ${\rm Hilb}^2S$ is obtained from the quotient $S \times S / \mathbb{Z}_2$ by blowing up the diagonal. Now if we combine the map $M = {\rm Hilb}^2S \to S \times S / \mathbb{Z}_2$ with the map $f : S \to \mathbb{CP}^1$ we see that we can map each point of $M$ to a pair of unordered points in $\mathbb{CP}^1$ which we can think of as the zeros of a quadratic polynomial, hence we get an induced map $\pi : M \to \mathbb{CP}^2$. This is in fact an example of holomorphic symplectic fibration.\\

To complete our example we should discuss what sort of singularities can be obtained and what the discriminant locus looks like. Matsushita \cite{mat} determines the types of singularities that can occur (except over a finite set of points) in a holomorphic symplectic fibration of the type we are considering. In fact most of the singularities listed in \cite{mat} are realised by taking the Hilbert scheme of two points on a K3 elliptic fibration in the previously described way. From this one should be able to the determine the monodromy representation on the fibre homology. It is also known the type of links that can occur as the result of placing a $3$-sphere around a singularity of a curve in $\mathbb{CP}^2$ \cite{boi}. For example a singularity of the form $z_1^p + z_2^q = 0$ (in affine coordinates) where $p$ and $q$ are coprime leads to a $(p,q)$-torus knot.

\section{Semi-flat coassociative fibrations}\label{secsemiflat}
We consider a relatively simple class of coassociative fibration we call semi-flat (perhaps more correctly $4/7$ flat) which extend the similar notion of a semi-flat special Lagrangian fibration.

\begin{defn}
Let $X$ be a $G_2$-manifold. We say that $X$ is {\em semi-flat} if $X$ is a coassociative fibration such that the fibres are flat tori.
\end{defn}
Since flat tori are compact Hyperk\"ahler manifolds, Proposition \ref{riemfibr} tells us that a semi-flat fibration $\pi : X \to B$ is also a Riemannian submersion. We call $V = {\rm ker}(\pi_*)$ the vertical distribution and the corresponding distribution of normals $V^\perp$ will be called the horizontal distribution.\\





\begin{lemp}\label{formofphi}
Let $X$ be a $G_2$-manifold. Suppose $V$ is a distribution of coassociative subspaces and $V^\perp$ the corresponding orthogonal distribution of associative subspaces. Given a local oriented orthonormal frame $\{e_j\}$, $1\le j \le 4$ for $V$ there exists an orthonormal frame $\{a_i\}$, $1\le i \le 3$ for $V^\perp$ such that the $3$-form $\phi$ and $4$-form $\psi$ have the following forms with respect to the corresponding coframe:
\begin{equation}\label{standard}
\begin{aligned}
\phi  &=  a^{123} + a^1 \! \! \wedge \! (e^{12}+e^{34}) + a^2 \! \! \wedge \! (e^{13} - e^{24}) + a^3 \! \! \wedge \! (-e^{14} - e^{23}) \\
\psi  &=  e^{1234} + a^{23} \! \! \wedge \! (e^{12}+e^{34}) + a^{31} \! \! \wedge \! (e^{13} - e^{24}) + a^{12} \! \! \wedge \! (-e^{14} - e^{23}).
\end{aligned}
\end{equation}
\begin{proof}
Follows since $G_2$ is transitive on the set of associative (or coassociative) subspaces with stabiliser ${\rm SO}(4)$ acting on the coassociative subspace by the standard representation \cite{harvlaw}.
\end{proof}
\end{lemp}

\begin{prop}\label{localaction}
Let $\pi : X \to B$ be a semi-flat $G_2$-manifold. Then locally $X$ has a $T^4$-action of $G_2$ isomorphisms. That is for any $b \in B$ there is a neighborhood $b \in U \subseteq B$ such that $\pi^{-1}(U)$ has a $T^4$-action preserving the $G_2$-structure and the orbits are the coassociative fibres.
\begin{proof}
Over a sufficiently small neighborhood in $B$ we can find a commuting frame $f_1, \dots , f_4$ for the vertical distribution $V$ such that for each fibre $F$ on which the frame is defined, flow along the integral curves of $f_1, \dots , f_4$ provides a diffeomorphism $\mathbb{R}^4/\mathbb{Z}^4 \to F$. This defines a local $T^4$-action on $X$ such that the orbits are fibres of $\pi$. The action clearly preserves the metric $g_V$ on the fibres.\\

Note that in general the frame $f_1 , \dots , f_4$ is not orthonormal. However we can make a change of frame to an orthonormal frame $e_1 , \dots , e_4$ such that the change of frame depends only on the base. Therefore $e_1 , \dots , e_4$ is a commuting orthonormal frame for $V$. Now let $a_1 , \dots , a_3$ be a frame for $V^\perp$ corresponding to the frame $e_1 , \dots , e_4$ for $V$ such that the associative $3$-form $\phi$ has the form in (\ref{standard}). For each fibre $F$ we find $\iota_{a_i} \phi |_F$ is closed (since the $e_i$ commute), so $a_i$ represents a deformation through coassociative fibres. But since the base exhausts all deformations it must be that the $a_i$ are lifts of a frame on $B$. It follows that the $a_i$ are invariant under the $T^4$-action. But it is clear that the $e_i$ are also $T^4$-invariant. Thus $\phi$ is invariant by (\ref{standard}).
\end{proof}
\end{prop}

\begin{prop}
Let $X$ be a semi-flat $G_2$-manifold. The horizontal distribution is integrable.
\begin{proof}
As in the proof of Proposition \ref{localaction} let $e_1, \dots , e_4$ be a commuting local frame for $V$ and $a_1, \dots , a_3$ the corresponding orthonormal frame for $V^\perp$ such that the $3$-form $\phi$ has the form in (\ref{standard}). Recall also that $a_1 , \dots , a_3$ is the horizontal lift of a frame on the base. If we let $e^1 , \dots , e^4$ and $a^1 , \dots , a^3$ be the corresponding coframes then the $a^i$ are the pull-back of a coframe on the base.\\

We now introduce a bi-grading on differential forms as follows. We identify $V^*$ and $(V^\perp)^*$ as subbundles of $T^*X$. Then a differential form of type $(p,q)$ is a section of $\wedge^p (V^\perp)^* \otimes \wedge^q V^*$. Since the $a^i$ the pull-back of forms on the base the $da^i$ are $(2,0)$-forms. Moreover since the $e_i$ are a commuting frame we have that $de^i$ is a sum of $(2,0)$ and $(1,1)$ terms.\\

It follows from (\ref{standard}) that $\phi$ is a sum of $(3,0)$ and $(1,2)$ terms and $\psi$ is a sum of $(2,2)$ and $(0,4)$ terms. Let us write $\psi = \psi_{(0,4)} + \psi_{(2,2)}$. We find that $d\psi_{(0,4)}$ is the sum of $(1,4)$ and $(2,3)$ terms while $d\psi_{(2,2)}$ is of type $(3,2)$. Therefore for $\psi$ to be closed we must have $d\psi_{(0,4)} = d\psi_{(2,2)} = 0$. In particular the $(2,3)$ component of $d(\psi_{(0,4)})$ vanishes. We find
\begin{equation}
\begin{aligned}
\left[d(\psi_{(0,4)}) \right]_{(2,3)} &= (de^1)_{(2,0)} \wedge e^{234} - e^1 \wedge (de^2)_{(2,0)} \wedge e^{34} \\
& \; \; \; + e^{12} \wedge (de^3)_{(2,0)} \wedge e^4 - e^{123} \wedge (de^4)_{(2,0)}.
\end{aligned}
\end{equation}
Setting this to zero implies that $(de^1)_{(2,0)} = \dots = (de^4)_{(2,0)} = 0$, hence the horizontal distribution is integrable.
\end{proof}
\end{prop}


Let $w \in U \subset \mathbb{R}^3$ be local coordinates for a leaf of the horizontal distribution. We can take $U$ sufficiently small that we have an embedded submanifold $i:U \to X$. Define a map $f:U \times T^4 \to X$ by $f(w,x) = x \cdot i(w)$ where $x$ acts on $i(w)$ by the local $T^4$-action. It is immediate that $f$ is an immersion sending $T(T^4)$ to the vertical distribution and $TU$ to the horizontal. Moreover by sufficiently restricting $U$ we may assume $f$ is injective. Let $\pi : X \to U$ denote the locally defined projection. The fibres of $\pi$ are coassociative submanifolds. Therefore we can identify $U$ as the local moduli space of coassociative deformations of fibres as in Proposition \ref{riemfibr}.\\

The metric $g$ has the form $g = g_V + g_{V^\perp}$ where $g_V$ and $g_{V^\perp}$ are metrics on the vertical and horizontal respectively. Now since $g$ is $T^4$-invariant we have that $g_V$ is a flat metric on each orbit and $g_{V^\perp}$ is the pull-back under $\pi$ of a metric $g_U$ on the base $U$ which is related to the $L^2$ moduli space metric by (\ref{basemetric}).


\section{Construction of semi-flat $G_2$-manifolds}\label{secconstruct}
We find an equivalent local characterisation of semi-flat $G_2$-manifolds in terms of minimal submanifolds.\\

Suppose for the moment that we have a semi-flat $G_2$-manifold. As usual let us take an invariant local frame $a_1,a_2,a_3,e_1, \dots , e_4$ such that $\phi$ has the form in (\ref{standard}). Let $b = (b^1,b^2,b^3)$ denote local horizontal coordinates for the base $B$ and let $x^i \in \mathbb{R}/\mathbb{Z}$, $1 \le i \le 4$ be standard coordinates for $T^4$. So the $db^\mu$ are $(1,0)$-forms and the $dx^i$ are $(0,1)$-forms. Hence we see that we can uniquely write $\phi$ as
\begin{equation}
\phi = dvol_B + db^\mu \wedge \theta_\mu
\end{equation}
where the $\theta_\mu$ are $(0,2)$-forms. It is also clear that $\theta_\mu|_{T_b}$ is a harmonic self-dual $2$-form on $T_b$ representing the deformation of $T_b$ in the $\partial / \partial b^\mu$ direction.\\

Thinking of $B$ as the moduli space of coassociative deformations we recall that there is a locally defined function $u : B \to H^2(T,\mathbb{R})$ such that $du$ is the $H^2(T,\mathbb{R})$-valued $1$-form $b \mapsto [ \theta_\mu |_{T_b} ]db^\mu$. We can give $H^2(T,\mathbb{R})$ coordinates $a_{ij}$ such that $a_{ij}$ corresponds to the cohomology class $[a_{ij}dx^{ij}]$. Therefore we have $6$ functions $u_{ij}$ on $B$ such that $u(b) = [u_{ij}(b)dx^{ij}]$. We then have
\begin{eqnarray*}
du &=& d[u_{ij}dx^{ij}]\\
&=& \left[ \frac{\partial u_{ij}}{\partial b^\mu} dx^{ij} \right] db^\mu \\
&=& [\theta_\mu |_{T_b}] db^\mu.
\end{eqnarray*}
Hence $[\theta_\mu |_{T_b}] = [\tfrac{\partial u_{ij}}{\partial b^\mu}(b) dx^{ij} ]$, in fact since the $\theta_\mu$ have no $db^\mu$ terms and are constant with respect to fibre coordinates (being harmonic) we define
\begin{equation}\label{theta}
\theta_\mu = \frac{\partial u_{ij}}{\partial b^\mu} dx^{ij}.
\end{equation}

By (\ref{basemetric}) the base metric $g_B$ and $L^2$ moduli space metric on $B$ are related by $g_{L^2} = 2{\rm vol}(T^4)g_B$. We can explicitly integrate the $L^2$ metric. If $g_B = g_{\mu \nu}db^\mu db^\nu$ then
\begin{eqnarray*}
g_{\mu \nu} &=& \dfrac{1}{2{\rm vol(T^4)}}\int_{T^4} u_* \left( \frac{\partial}{\partial b^\mu} \right) \wedge u_* \left( \frac{\partial}{\partial b^\nu} \right) \nonumber \\
&=& \dfrac{1}{2{\rm vol(T^4)}}\int_{T^4} \frac{\partial u_{ij}}{\partial b^\mu}dx^{ij} \wedge \frac{\partial u_{kl}}{\partial b^\nu}dx^{kl} \nonumber \\
&=& \dfrac{1}{2{\rm vol(T^4)}}\epsilon^{ijkl} \frac{\partial u_{ij}}{\partial b^\mu} \frac{\partial u_{kl}}{\partial b^\nu}
\end{eqnarray*}
where $\epsilon^{abcd}$ is antisymmetric and $\epsilon^{1234}=1$. So we have
\begin{equation}\label{metric}
2g_B(A,B){\rm vol(T^4)}dx^{1234} = u_*(A) \wedge u_*(B),
\end{equation}
so $g_B$ is essentially the pull-back of the wedge product.\\

We are now in a position to reverse the construction. Suppose $B$ is an oriented $3$-manifold with a function $u : B \to H^2(T,\mathbb{R})$, $u = [u_{ij}dx^{ij}]$. We assume that $u_*$ sends the tangent spaces of $B$ into maximal positive definite subspaces of $H^2(T,\mathbb{R})$. Choose a positive constant $\tau$ representing the volume of the coassociative fibres. We may pull back the intersection form to define a positive definite metric $h$ on $B$ given by 
\begin{equation}\label{metrich}
2h(A,B)\tau dx^{1234} = u_*(A) \wedge u_*(B).
\end{equation}
Let $dvol_h$ denote the volume form for this metric. We may then define self-dual $2$-forms $\theta_\mu$ by equation (\ref{theta}), hence we can define a $3$-form $\phi$ on $B \times T^4$ by
\begin{equation}\label{constructedphi}
\phi = dvol_h + db^\mu \wedge \theta_\mu.
\end{equation}
We can easily verify that $\phi$ is closed by noting that from (\ref{theta}) we have $(\partial/\partial b^\nu)\theta_\mu = (\partial/\partial b^\mu)\theta_\nu$.\\

It is clear that for any $u$ and $\tau$ as above, $\phi$ has the correct algebraic form for an associative $3$-form. By (\ref{gdvol1}), $\phi$ determines a metric $g$ on $X$ and a corresponding volume form $dvol_X = dvol_B \wedge dvol_T$. It follows from (\ref{constructedphi}) that the induced metric $g$ agrees with $h$ on the horizontal distribution so that $dvol_B = \, dvol_h$. Moreover one can further show $dvol_T = \tau dx^{1234}$.\\

Now consider $\psi = *\phi$. Since the $\theta_\mu$ are self-dual $2$-forms on each fibre we find
\begin{equation}
\psi = dvol_T + *_3 db^\mu \wedge \theta_\mu
\end{equation}
where $*_3$ denotes the Hodge star with respect to $g$ restricted to the base. We can see that a necessary condition for $\psi$ to be closed is that $\tau$ is constant. Indeed if $\psi$ is closed then since it is a calibrating form and the fibres of $X$ are isotopic calibrated submanifolds they must have equal volume. Therefore assume $\tau$ is constant. We calculate
\begin{eqnarray*}
d\psi &=& 0 + d(*_3 db^\mu) \wedge \theta_\mu + *_3 db^\mu \wedge d\theta_\mu \\
&=& \Delta b^\mu dvol_B \wedge \theta_\mu + g^{\mu \nu} \iota_\nu dvol_B \wedge db^\gamma \wedge \frac{\partial^2 u_{ij}}{\partial b^\gamma \partial b^\mu} dx^{ij} \\
&=& \Delta b^\mu dvol_B \wedge \frac{\partial u_{ij}}{\partial b^\mu} dx^{ij} + g^{\mu \nu} dvol_B \wedge \frac{\partial^2 u_{ij}}{\partial b^\nu \partial b^\mu} dx^{ij} \\
&=& \left( \Delta b^\mu \frac{\partial u_{ij}}{\partial b^\mu} + g^{\mu \nu} \frac{\partial^2 u_{ij}}{\partial b^\mu \partial b^\nu} \right) dvol_B \wedge dx^{ij}.
\end{eqnarray*}
Hence $\psi$ is closed if and only if for each $i,j$ we have
\begin{equation}\label{uequ}
g^{\mu \nu} \frac{\partial^2 u_{ij}}{\partial b^\mu \partial b^\nu} + \Delta b^\mu \frac{\partial u_{ij}}{\partial b^\mu} = 0.
\end{equation}
Note that the Laplacian on functions on the base $\Delta = *^{-1}d*d = -\delta d$ is given by
\begin{equation}
\Delta f = g^{\mu \nu} \left(\frac{\partial^2 f}{\partial b^\mu \partial b^\nu} - {\Gamma^\gamma}_{\mu \nu} \frac{\partial f}{\partial b^\gamma} \right).
\end{equation}
Where ${\Gamma^\gamma}_{\mu \nu}$ are the Christoffel symbols
\begin{equation}
{\Gamma^\gamma}_{\mu \nu} = \frac{1}{2}g^{\gamma \sigma}( \partial_\mu g_{\nu \sigma} + \partial_\nu g_{\mu \sigma} - \partial_\sigma g_{\mu \nu} ).
\end{equation}
In particular, applied to a coordinate function we have $\Delta b^\gamma = -g^{\mu \nu}{\Gamma^\gamma}_{\mu \nu}$. Substituting this into equation (\ref{uequ}) we get
\begin{equation}
\Delta u_{ij} = 0.
\end{equation}
This says that the map $u : B \to H^2(T,\mathbb{R})$ is harmonic where $B$ is given the metric $g$. However the pull-back metric on $B$ induced by $u$ differs from $g$ only by a constant so equivalently $u$ is harmonic with respect to the induced metric. Another way of saying this is that the map $u$ is a minimal immersion or that $B$ is locally embedded as a minimal $3$-submanifold \cite{eells}. 
\begin{thm}\label{thethm}
Let $B$ be an oriented $3$-manifold and $u : B \to \wedge^2 \mathbb{R}^4$ a map with the property that $u$ maps the tangent spaces of $B$ into maximal positive definite subspaces of $\wedge^2 \mathbb{R}^4$ and let $\tau$ be a positive constant. Let $h$ be the pull-back metric defined in equation (\ref{metrich}) with volume $dvol_h$. Let $X = B \times \left( \mathbb{R}/\mathbb{Z} \right)^4$ and define $\phi \in \Omega^3(X,\mathbb{R})$ by
\begin{equation}
\phi = \, dvol_h + du,
\end{equation}
where $u$ is considered as a $2$-form on $X$. Then $(X,\phi)$ is a semi-flat $G_2$-manifold if and only if $u$ is a minimal immersion. Moreover every semi-flat $G_2$-manifold locally has this form.
\end{thm}

We can improve on Theorem \ref{thethm} by determining the global properties of semi-flat $G_2$-manifolds. For completeness we may define a {\em locally semi-flat} $G_2$-manifold as a locally trivial $T^4$-fibration $\pi : X \to M$ with semi-flat local trivialisations. We leave it to the reader fill in the details which ammounts to incorporating monomdromy representations.

\section{Cylindrical semi-flat $G_2$-manifolds}\label{seccylsem}
We will show that the semi-flat $G_2$ equations reduce to the Monge-Amp\`ere equation if we assume that the resulting $G_2$-manifold is cylindrical. This coincides with the result of Hitchin \cite{hit} in which special Lagrangian fibrations with flat fibres are produced from the Monge-Amp\`ere equation.\\

Let $X$ be a semi-flat $G_2$-manifold over a base $B$. Suppose that $X = Y \times T$ as a Riemannian manifold where $T$ is a circle and $Y$ is a Calabi-Yau manifold. Moreover we suppose that the circle factor generates a $T^1$ subgroup of the $T^4$-action on $X$. Then $Y$ is a semi-flat special Lagrangian fibration over the same base. In this situation we say that $X$ is {\em cylindrical}.\\

Let us write the $4$-torus $T^4$ as a product $T^4 = T^3 \times T^1$. Using duality between $H^1(T^3,\mathbb{R})$ and $H^2(T^3,\mathbb{R})$ we equip the space $H^1(T^3,\mathbb{R}) \oplus H^2(T^3,\mathbb{R})$ with a split signature inner product $\langle \, , \, \rangle$ and symplectic form $\omega$. For $\zeta^1 , \eta^1 \in H^1(T^3,\mathbb{R})$, $\zeta^2, \eta^2 \in H^2(T^3,\mathbb{R})$ we have
\begin{equation}
\begin{aligned}
\langle (\zeta^1 , \zeta^2) , (\eta^1 , \eta^2) \rangle =& \, \frac{1}{2}(\zeta^1(\eta^2) + \eta^1(\zeta^2))\\
\omega ( (\zeta^1 , \zeta^2) , (\eta^1 , \eta^2) ) =& \, \frac{1}{2}(\zeta^1(\eta^2) - \eta^1(\zeta^2)).
\end{aligned}
\end{equation}

We have an isometry $e : H^1(T^3,\mathbb{R}) \oplus H^2(T^3,\mathbb{R}) \to H^2(T^4,\mathbb{R})$ as follows:
\begin{equation}
e(\alpha, \beta) = (\alpha \smallsmile [dt] + \beta)/\sqrt{2}
\end{equation}
where $t \in \mathbb{R}/\mathbb{Z}$ is the standard coordinate for $T$.\\

Given a cylindrical semi-flat $G_2$ manifold $X = Y \times T$ we have the corresponding minimal immersion $\tilde{u} : B \to H^2(T^4 , \mathbb{R})$. We can express $\tilde{u}$ as the composition $\tilde{u} = \sqrt{2}e( u,v)$ where $u : B \to H^1(T^3,\mathbb{R})$, $v: B \to H^2(T^3,\mathbb{R})$. From \cite{hit} we know that $(u,v)$ locally embeds $B$ as a Lagrangian submanifold. 

Let $x^1,x^2,x^3$ denote coordinates on $T^3$ and $t = x^0$ a coordinate for $T$. The maps $u,v$ are local diffeomorphisms hence writing $u = [u_i dx^i]$ and $v = [\tfrac{1}{2!} \epsilon_{ijk}v^idx^{jk}]$ we may take either the $u_i$ or $v^j$ as coordinates on $B$. Moreover since $B$ is Lagrangian it is locally the graph of a gradient, that is there exist functions $\phi,\psi$ on $B$ such that $v^i = \frac{\partial \phi}{\partial u_i}$, $u_i = \frac{\partial \psi}{\partial v^i}$. Now we write $\tilde{u} = [\tilde{u}_{ij}dx^{ij}]$ so that
\begin{equation}
\begin{aligned}
\tilde{u}_{i0} &= u_i, \\
\tilde{u}_{ij} &= \tfrac{1}{2}\epsilon_{ijk}v^k, \; \; i,j,k > 0.
\end{aligned}
\end{equation}
If we check the formula for the metric induced by $\tilde{u}$ we find that up to a multiple it is given by
\begin{equation}
g = \frac{\partial^2 \phi}{\partial u_i \partial u_j} du_i du_j = \frac{\partial^2 \psi}{\partial v^i \partial v^j} dv^i dv^j.
\end{equation}
Now starting with the relation $v^j = \frac{\partial \phi}{\partial u_j}$ we find $dv^j = \frac{\partial^2 \phi}{\partial u_j \partial u_k}du_k$. Substituting this into the expression for $g$ we find that $h^{ij} = \frac{\partial^2 \psi}{\partial v^i \partial v^j}$ is the inverse of $g_{ij} = \frac{\partial^2 \phi}{\partial u_i \partial u_j}$, that is $g_{ij}h^{jk} = \delta^k_i$. Let us introduce some notation: $\partial_i = \tfrac{\partial}{\partial u_i}$, $\partial^j = \tfrac{\partial}{\partial v^j}$, $\phi_{ij \dots k} = \partial_i \partial_j \dots \partial_k \phi$, $\psi^{ij \dots k} = \partial^i \partial^j \dots \partial^k \psi$. We also note that $\partial^j = \psi^{jk}\partial_k$. Now we can calculate the Christoffel symbols in the $u_i$ coordinates:
\begin{equation}
{\Gamma^k}_{ij} = \frac{1}{2}\psi^{km}\psi_{mij}.
\end{equation}
Now we calculate $\Delta \tilde{u}_{ab}$. First suppose $b=0$ so that $\tilde{u}_{ab} = u_a$. Then we find
\begin{eqnarray}
\Delta u_a &=& \psi^{ij}\left( 0 - \frac{1}{2}\psi^{km}\phi_{mij}\delta^a_k \right) \nonumber \\
&=& -\frac{1}{2}\psi^{ij}\psi^{am}\phi_{mij}. \label{laplacian}
\end{eqnarray}
If we take the relation $\phi_{ij}\psi^{jk} = \delta^k_i$ and differentiate we find
\begin{equation}
\phi_{mij}\psi^{jk} = -\phi_{ij}\phi_{mr} \psi^{rjk}.
\end{equation}
Substituting into equation (\ref{laplacian}) we find that
\begin{equation}
\Delta u_a = \frac{1}{2}\phi_{ij}\psi^{aij}.
\end{equation}
Similarly we find
\begin{equation}
\Delta v^a = \frac{1}{2}\psi^{ij}\phi_{aij}.
\end{equation}
Therefore the $G_2$ equations in this case reduce to
\begin{equation}
\begin{aligned}
\phi_{ij}\psi^{aij} =& 0, \\
\psi^{ij}\phi_{aij} =& 0.
\end{aligned}
\end{equation}

Now let us recall Jacobi's formula in the case where the matrix valued function $\phi$ is invertible:
\begin{equation*}
d {\rm det}(\phi) = {\rm det}(\phi) {\rm tr}(\phi^{-1}d\phi).
\end{equation*}
Therefore the Monge-Amp\`ere equation $d {\rm det}(\phi) = 0$ is equivalent to ${\rm tr}(\phi^{-1}d\phi) = 0$. But
\begin{eqnarray*}
{\rm tr}(\phi^{-1}d\phi) &=& {\rm tr}( \psi d\phi) \\
&=& \psi^{ij}\phi_{aij} du_a.
\end{eqnarray*}
Similarly we can interchange the roles of $\phi$ and $\psi$. Hence the $G_2$ equations in this case are equivalent to the Monge-Amp\`ere equation. 
\begin{rem}
We note that in \cite{hit} the Monge-Amp\`ere equation is also shown to be equivalent to $B$ being calibrated with respect a calibrating form that is a linear combination of the volume forms of $H^1(T^3,\mathbb{R})$ and $H^2(T^3,\mathbb{R})$. This agrees with the fact that $B$ is minimally immersed.
\end{rem}

\section{Reduction to surface equations}\label{secredtosurf}
We impose an additional degree of symmetry on a semi-flat $G_2$-manifold. The additional symmetry is shown to correspond to the base locally having the structure of a minimal cone in $\mathbb{R}^{3,3}$ which in turn is equivalent to a minimal surface in the quadric of unit vectors.\\

Let $\pi:X \to B$ be a semi-flat $G_2$-manifold constructed from a minimal immersion $u: B \to \mathbb{R}^{3,3}$. We suppose there is a vector field $U$ on $X$ such that $U$ commutes with the $T^4$-action. It is not possible for such a vector field to preserve the $3$-form $\phi$ up to scale, that is $\mathcal{L}_U \phi = \lambda \phi$ for some non-vanishing function $\lambda$, for in this case the fibres of the $T^4$-fibration would not have constant volume. Therefore we consider a slightly different symmetry. We suppose that 
\begin{equation}\label{scale}
\mathcal{L}_U\phi = \lambda \phi + 2 \lambda dvol_B
\end{equation}
as such a symmetry will preserve the volume of the fibres.\\

Let $U =V+W = V^\mu \tfrac{\partial}{\partial b^\mu} + W^i \tfrac{\partial}{\partial x^i} $. For $U$ to commute with the $T^4$-action we must have the $V^\mu$ and $W^i$ are independent of $x$. We will show that $W$ is a vector field generated by the $T^4$-action, hence we need only consider $V$. 

Recall that locally a semi-flat $G_2$-manifold $X = B \times T^4$ with coordinates $(b,x)$ has the $3$-form
\begin{equation}
\phi = dvol_B + du_{ij} \wedge dx^{ij} = dvol_B + db^\mu \wedge \theta_\mu
\end{equation}
where $u : B \to H^2(T^4,\mathbb{R}/\mathbb{Z})$ is a minimal immersion and $\theta_\mu = \tfrac{\partial u_{ij}}{\partial b^\mu}dx^{ij}$. For simplicity we will take ${\rm vol}(T^4)=1$.\\

Since $\phi$ is closed, the condition on $U$ is that $\mathcal{L}_U\phi = d(\iota_U\phi)= \lambda \phi + 2\lambda dvol_B$. We find that 
\begin{equation*}
\iota_U\phi = \iota_V dvol_B + V^\mu \theta_\mu - db^\mu \wedge \iota_W \theta_\mu
\end{equation*}
and that
\begin{equation*}
d(\iota_U \phi) = {\rm div}(V)dvol_B + d(V^\mu\theta_\mu) + 2 db^\mu \wedge d(W^i \frac{\partial u_{ij}}{\partial b^\mu}dx^j ).
\end{equation*}
Equating this to (\ref{scale}) we find
\begin{equation}
\begin{aligned}
{\rm div}(V) &= 3\lambda, \\
d(V^\mu \theta_\mu) &= \lambda db^\mu \wedge \theta_\mu, \\
\frac{\partial W^i}{\partial b^\mu}\frac{\partial u_{ij}}{\partial b^\nu} &= \frac{\partial W^i}{\partial b^\nu}\frac{\partial u_{ij}}{\partial b^\mu}.
\end{aligned}
\end{equation}
We can show that $W$ must be independent of the base variables. This follows from the algebraic fact that if $\theta_1,\theta_2,\theta_3$ are a basis of self-dual $2$-forms and $A_1,A_2,A_3$ are vectors such that $i_{A_i}\theta_j = i_{A_j}\theta_i$ then $A_i = 0$. Hence $W$ is a vector field coming from the $ T^4$-action. Therefore we ignore $W$.\\

Now the equations for $V$ are $3\lambda = {\rm div}(V)$ and $d(V^\mu \theta_\mu) = \lambda db^\mu \wedge \theta_\mu = \lambda du$. Taking exterior derivatives we find $0 = d\lambda \wedge du$ hence $\tfrac{\partial \lambda}{\partial b^\mu}\tfrac{\partial u}{\partial b^\nu} = \tfrac{\partial \lambda}{\partial b^\nu}\tfrac{\partial u}{\partial b^\mu}$, that is $\tfrac{\partial \lambda}{\partial b^\mu}\theta_\nu = \tfrac{\partial \lambda}{\partial b^\nu}\theta_\mu$. But $\{\theta_\mu\}$ are linearly independent so we have $d\lambda = 0$. Therefore the second equation for $V$ becomes $d(V^\mu\theta_\mu - \lambda u)=0$ or $u_*(V) = \lambda u + c$ where $c$ is constant. There are now two cases to consider; when $\lambda =0$ and $\lambda \ne 0$.

The $\lambda = 0$ case can readily be shown to correspond to minimal surfaces in $\mathbb{R}^{2,3}$. Such minimal surfaces correspond locally to holomorphic maps $\tau : \Sigma \to Q$ from a Riemann surface into an open subset of a quadric given by $Q = \{ v \in \mathbb{C} \otimes \mathbb{R}^{2,3} \; | \; \langle v , v \rangle = 0, \; \langle v , \overline{v} \rangle > 0 \}$. There is a local Weierstrass representation for the corresponding minimal immersion $\phi$ \cite{li}:
\begin{equation*}
\phi(z) = \phi(0) + {\rm Re}\int_0^z \tau(\zeta) d\zeta.
\end{equation*}
Now assume $\lambda \ne 0$. Then we can redefine $u$ to absorb the constant $c$ so we have $u_*(V) = \lambda u$. Now we can rescale $V$ such that $u_*(V) = u$. If the vector field $V$ vanishes at a point $b \in B$ then $u(b) = u_*(V_b)=0$, but the map $u:B \to H^2(T^4,\mathbb{R})$ is an immersion so $V$ vanishes at isolated points. Away from the zeros of $V$ we may find local coordinates $(x,y,t) \in \Sigma \times I$ such that $V = \tfrac{\partial}{\partial t}$. Hence we have $\tfrac{\partial u}{\partial t}(x,y,t) = u(x,y,t)$. The solution is of the form $u(x,y,t) = u(x,y)e^t$. The induced metric $g$ on $B$ has the property that $g(x,y,t) = e^{2t}g(x,y)$, hence $\tfrac{\partial}{\partial t}$ satisfies ${\rm div}(\tfrac{\partial}{\partial t})=3$ as required.\\

We will attempt to find local coordinates that diagonalise the metric on $B$. By changing the local slice $\Sigma \to B$ along the $t$ direction we have freedom $u(x,y) \mapsto u(x,y)e^{\rho(x,y)}$ where $\rho$ is an arbitrary smooth function on $B$. We calculate (on $t=0$)
\begin{equation*}
2 g \left( \frac{\partial}{\partial t} , \frac{\partial}{\partial x} \right) dx^{1234} = e^{2\rho}( u \wedge \frac{\partial u}{\partial x} + u \wedge u \frac{\partial \rho}{\partial x} )
\end{equation*}
and similarly for $g \left( \frac{\partial}{\partial t} , \frac{\partial}{\partial y} \right)$. Let us define functions $r,s$ by
\begin{equation}
\begin{aligned}
r u \wedge u &= -u \wedge \frac{\partial u}{\partial x} \\
s u \wedge u &= -u \wedge \frac{\partial u}{\partial y}.
\end{aligned}
\end{equation}
Note that this is possible because $2g(\tfrac{\partial}{\partial t} , \tfrac{\partial}{\partial t} )dx^{1234} = u \wedge u \neq 0$. Then we can locally find a function $\rho(x,y)$ such that $g(\tfrac{\partial}{\partial t} , \tfrac{\partial}{\partial x} ) = g(\tfrac{\partial}{\partial t} , \tfrac{\partial}{\partial y}) = 0$ if and only if $\tfrac{\partial r}{\partial y} = \tfrac{\partial s}{\partial x}$. This follows easily from the definitions of $r$ and $s$. Therefore our metric now has the form
\begin{equation*}
g(x,y,t) = e^{2t}(c(x,y)dt^2 + g_{\Sigma}(x,y) ).
\end{equation*}
Setting $r = e^t$ we may write this as
\begin{equation*}
g(x,y,r) = c(x,y)dr^2 + r^2g_{\Sigma}(x,y).
\end{equation*}
Let us also note that $2c dx^{1234} = u \wedge u$ so that $2\tfrac{\partial c}{\partial x} dx^{1234} =  2 u \wedge \tfrac{\partial u}{\partial x} = 0$ and similarly for $y$. Thus $c$ is constant and by scaling we can assume $c=1$. Therefore the metric is
\begin{equation}
g = dr^2 + r^2 g_{\Sigma}
\end{equation}
and our minimal $3$-fold is a cone.\\

For a manifold $M$ with (possibly indefinite) metric $g$ let $(CM,\hat{g})$ denote the cone where $CM = M \times (0,\infty)$ , $\hat{g} = dr^2 + r^2g$. Given a map $\phi : (M,g) \to (N,h)$ we define the {\em radial extension} $\hat{\phi} : CM \to CN$ by $\hat{\phi}(x,r) = (\phi(x),r)$.

\begin{prop}
The radial extension $\hat{\phi} : CM \to CN$ is minimal if and only if $\phi : M \to N$ is minimal.
\begin{proof}
First we note that $\phi$ is a Riemannian immersion if and only if $\hat{\phi}$ is a Riemannian immersion. Let us use coordinates $x^1, \dots , x^m$ on $M$ and let $r = x^0$. We use the convention that indices $i,j,k,\dots $ do not take the value $0$. Likewise give $CN$ coordinates $r=y^0,y^1, \dots , y^n$. We have
\begin{equation}
\begin{aligned}
\hat{g}_{00} &= 1, \; \; \; & \hat{g}_{0i} &= 0, \; \; \; & \hat{g}_{ij} &= r^2 g_{ij} \\
\hat{g}^{00} &= 1, \; \; \; & \hat{g}^{0i} &= 0, \; \; \; & \hat{g}^{ij} &= \frac{1}{r^2} g^{ij}.
\end{aligned}
\end{equation}
We readily verify the following relation between the Christoffel symbols on $M$ and $CM$:
\begin{equation}
\begin{aligned}
\leftexp{CM}{{\Gamma^k}_{ij}} &= \leftexp{M}{{\Gamma^k}_{ij}}, \; \; \; &
\leftexp{CM}{{\Gamma^0}_{ij}} &= -rg_{ij}, \; \; \; &
\leftexp{CM}{{\Gamma^k}_{0j}} &= \frac{1}{r}\delta^k_j, \\
\leftexp{CM}{{\Gamma^k}_{00}} &= 0, \; \; \; &
\leftexp{CM}{{\Gamma^0}_{0j}} &= 0, \; \; \; &
\leftexp{CM}{{\Gamma^0}_{00}} &= 0.
\end{aligned}
\end{equation}
There are similar relations for the Christoffel symbols on $CN$. The map $\hat{\phi}$ has the properties
\begin{eqnarray*}
\phi^0 &=& x^0 \\
\partial_0 \phi^\gamma &=& 0.
\end{eqnarray*}
Recall \cite{eells} that the tension field $\tau(\phi)$ for a map $\phi :M \to N$ is the section of $\phi^*(TN)$ obtained by taking the trace of the second fundamental form:
\begin{equation}
\tau^\gamma(\phi) = g^{ij}\left( \partial^2_{ij} \phi^\gamma - \leftexp{M}{{\Gamma^k}_{ij}} \partial_k \phi^\gamma + \leftexp{N}{{\Gamma^\gamma}_{\alpha \beta}} \partial_i \phi^\alpha \partial_j \phi^\beta \right).
\end{equation}
The map $\phi$ is harmonic if and only if $\tau(\phi)=0$. Similarly we have a torsion field $\tau(\hat{\phi})$. We calculate
\begin{eqnarray}
\tau^\gamma(\hat{\phi}) &=& \frac{1}{r^2} \tau^\gamma(\phi) \\
\tau^0(\hat{\phi}) &=& \frac{1}{r} g^{ij} \left( g_{ij} - h_{\alpha \beta}\partial_i \phi^\alpha \partial_j \phi^\beta \right).
\end{eqnarray}
The result follows.
\end{proof}
\end{prop} 

Let $\mathbb{R}^{p,q}$ denote $\mathbb{R}^n$ with a signature $(p,q)$ inner product. We say that a submanifold $X$ of $\mathbb{R}^{p,q}$ is a cone if $X$ is diffeomorphic to $\Sigma \times (0,\infty)$ such that $i(x,r) = ri(x,1)$ where $i$ is the inclusion $i : X \to \mathbb{R}^{p,q}$ and the induced metric on $X$ is of the form $dr^2 + r^2g_\Sigma$ where $g_\Sigma$ is independent of $r$. 

Restricting to $r=1$ we have an inclusion $ i : \Sigma \to Q \subset \mathbb{R}^{p,q}$ where $Q = \{ v \in \mathbb{R}^{p,q} | \; \langle v , v \rangle = 1 \} = {\rm O}(p,q)/{\rm O}(p-1,q)$. Conversely such a map defines a cone in $\mathbb{R}^{p,q}$. We thus have
\begin{corp}
There is a bijection between minimal cones with definite induced metric in $\mathbb{R}^{p,q}$ and minimal submanifolds of the quadric $Q = {\rm O}(p,q)/{\rm O}(p-1,q)$ with definite induced metric.
\end{corp}

In our situation we are looking for minimal surfaces $u : \Sigma \to {\rm O}(3,3)/{\rm O}(2,3)$.

\section{Minimal surfaces of signature $(3,3)$}\label{secquadrics2}

We will develop some general theory of minimal surfaces in a quadric. For the most part this is a straightforward generalisation of the definite signature case \cite{bol}.\\

Let $\mathbb{R}^{p,q}$ denote the $p+q$-dimensional vector space with bilinear form $\langle \; , \; \rangle$ of signature $(p,q)$ and let $Q_{h_0} = \{x \in \mathbb{R}^{p,q} \,| \,\langle x , x \rangle = h_0 \}$ where $h_0 = \pm 1$. Let $\Sigma$ be a connected oriented surface. An immersion $\phi : \Sigma \to Q_{h_0} \subset \mathbb{R}^{p,q}$ that sends each tangent space of $\Sigma$ to a positive definite subspace of $\mathbb{R}^{p,q}$ induces a metric and hence a conformal structure on $\Sigma$. These give $\Sigma$ the structure of a Riemann surface with compatible metric. The map $\phi$ is then a minimal immersion if and only if it is harmonic.

Let us regard $\phi$ as a vector valued function on $\Sigma$ such that $\langle \phi , \phi \rangle = h_0$. Let $D$ denote the trivial connection on $\mathbb{R}^{p,q}$ and $\nabla$ the induced Levi-Civita connection on $Q$. Let $z$ be a local holomorphic coordinate on $\Sigma$. We use $z$ and $\overline{z}$ subscripts to denote partial differentiation with respect to $\tfrac{\partial}{\partial z}$ and $\tfrac{\partial}{\partial \overline{z}}$. The harmonic equation for $\phi$ is
\begin{equation}
\nabla_{\overline{z}}{\phi_z} = 0.
\end{equation}
However we also have that $D = \nabla + \Pi$ where $\Pi$ is the second fundamental form which is valued in the normal bundle of $Q$, since $Q$ is given the induced metric. Therefore the harmonic equation for $\phi$ reduces to
\begin{equation}
\phi_{z\overline{z}} = \lambda \phi
\end{equation}
for some function $\lambda$. In fact we can show $\lambda = -\langle \phi_z , \phi_{\overline{z}} \rangle/h_0$. Note that since $\langle \phi , \phi \rangle = h_0$ we have $\langle \phi , \phi_z \rangle = \langle \phi , \phi_{\overline{z}} \rangle = 0$ and since $\phi$ is a Riemannian immersion $\langle \phi_z , \phi_z \rangle = \langle \phi_{\overline{z}} , \phi_{\overline{z}} \rangle = 0$ and $\langle \phi_z , \phi_{\overline{z}} \rangle = h_1$ defines the induced metric on $\Sigma$.\\

Following the definite signature case \cite{bol} we will inductively construct a sequence $\phi_0,\phi_1,\phi_2, \dots $. We define $\phi_0 = \phi$, $\phi_1 = \phi_z$ and we define $\phi_{i+1}$ from $\phi_i$ under the assumption that $h_i = \langle \phi_i , \overline{\phi_i} \rangle$ is non-vanishing by
\begin{equation}
\phi_{i+1} = (\phi_i)_z - \frac{\langle (\phi_i)_z , \overline{\phi_i} \rangle}{h_i}\phi_i.
\end{equation}
So $\phi_{i+1}$ is the projection of $(\phi_i)_z$ onto the orthogonal complement of $\phi_i$ with respect to the Hermitian form $h(X,Y) = \langle X , \overline{Y} \rangle$. Note that unlike the definite signature case it may be that $h_i$ vanishes even if $\phi_i \ne 0$.\\

Let $r \ge 1$ be an integer. We say $\phi$ has {\em isotropy of order $r$} if we can construct the sequence $\phi_0,\phi_1, \dots , \phi_r, \phi_{r+1}$ (so $h_1,\dots, h_r$ are non-vanishing) and $\langle \phi_i , \phi_i \rangle = 0$ for $1\le i \le r$. Since the $h_i$ are non-vanishing we may write $h_i = \epsilon_i H_i$ where $H_i = e^{2w_i}$ and $\epsilon_i = \pm 1$. We then have
\begin{prop}
Let $\phi$ be a minimal surface with isotropy of order $r \ge 1$. We have the following
\begin{itemize}
\item{The sequence $\overline{\phi_r}, \dots , \overline{\phi_1} , \phi_0 , \phi_1, \dots , \phi_r$ is pairwise orthogonal with respect to the Hermitian form,}
\item{$\phi_{r+1}$ and $\overline{\phi_{r+1}}$ are orthogonal to $\overline{\phi_r}, \dots , \overline{\phi_1} , \phi_0 , \phi_1, \dots , \phi_r$,}
\item{$\langle \phi_{r+1} , \phi_{r+1} \rangle$ defines a holomorphic degree $2r+2$-differential on $\Sigma$,}
\item{$(\phi_i)_z = \phi_{i+1} + 2(w_i)_z \phi_i$, for $1 \le i \le r$,}
\item{$(\phi_i)_{\overline{z}} = -(h_{i}/h_{i-1}) \phi_{i-1}$, for $1 \le i \le r+1$,}
\item{$2(w_i)_{z\overline{z}} = h_{i+1}/h_i - h_i/h_{i-1}$, for $1 \le i \le r$.}
\end{itemize}
\begin{proof}
These are all reasonably straightforward verifications nearly identical to the definite signature case.
\end{proof}
\end{prop}

\begin{defn}
Let $\phi:\Sigma \to Q_{h_0} \subset \mathbb{R}^{p,q}$ be a minimal surface. Suppose $\phi$ has isotropy of order $r$ and that $p+q = 2r+2$ or $2r+3$. We say $\phi$ is {\em superminimal} if $h_{r+1}=0$ and {\em superconformal} if $h_{r+1} \ne 0$.
\end{defn}

We now consider a superconformal or superminimal surface $\phi: \Sigma \to Q_{h_0} \subset \mathbb{R}^{p,q}$ in the case $p+q$ is even. So $\phi$ has isotropy of order $r$ such that $p+q = 2r+2$. Now $\overline{\phi_r}, \cdots , \overline{\phi_1} , \phi_0 , \phi_1, \cdots , \phi_r$ span a real codimension $1$ subspace which is not null, therefore there is a real vector valued function $\tilde{\phi}$ with $\langle \tilde{\phi} , \tilde{\phi} \rangle = \epsilon = \pm 1$ and complex valued function $q$ such that $\phi_{r+1} = q\tilde{\phi}$ and $\overline{\phi_{r+1}} = \overline{q}\tilde{\phi}$. Thus $h_{r+1} = \epsilon q \overline{q}$ and $\langle \phi_{r+1} , \phi_{r+1} \rangle = \epsilon q^2$ is holomorphic. Hence $q(dz)^{r+1}$ is a holomorphic $(r+1)$-differential. We then have the following set of equations:
\begin{equation} \label{signedtoda1}
\begin{aligned}
2(w_i)_{z\overline{z}} &= h_{i+1}/h_i - h_i/h_{i-1}, \; \; \; 1 \le i \le r-1, \\
2(w_r)_{z\overline{z}} &= \epsilon q\overline{q}/h_r - h_r/h_{r-1}.
\end{aligned}
\end{equation}
Notice that the superminimal case where $\phi_{r+1}=0$ corresponds to setting $q=0$. We can also write the equations directly in terms of the $w_i$. Let $\mu_i = \epsilon_i\epsilon_{i+1}$ for $1\le i \le r$ and $\mu_{r+1} = \epsilon \epsilon_{r}$. Then
\begin{equation}
\begin{aligned}
2(w_i)_{z\overline{z}} &= \mu_{i+1}e^{2w_{i+1} - 2w_i} -\mu_i e^{2w_i - 2w_{i-1}}, \; \; \; 1 \le i \le r-1, \\
2(w_r)_{z\overline{z}} &= \mu_{r+1}q\overline{q}e^{-2w_r} -\mu_r e^{2w_r - 2w_{r-1}}.
\end{aligned}
\end{equation}
We find this is a set of affine Toda equations for the affine Dynkin diagram $D^{(2)}_{r+1}$. The different choices of signs $\mu_i = \pm 1$ determine different real forms of the Toda equations.\\

Let us apply this to the case of a minimal surface $\phi$ into $Q = \{ x \in \mathbb{R}^{3,3} \, | \, \langle x , x \rangle = 1 \}$. In this case the maximum possible isotropy order is $2$ and we have $\mu_1 = 1$, $\mu_2 = -1$, $\mu_3 = 1$, $q$ is a holomorphic cubic differential and the equations become
\begin{equation}\label{d23toda}
\begin{aligned}
2(w_1)_{z\overline{z}} &= -e^{2w_2-2w_1} - e^{2w_1}, \\
2(w_2)_{z\overline{z}} &= q\overline{q}e^{-2w_2} + e^{2w_2 - 2w_1}.
\end{aligned}
\end{equation}
Let $a_1 = H_1 = e^{2w_1}$ and $a_2 = H_2/H_1 = e^{2w_2-2w_1}$. Then $a_1$ and $a_2$ are positive $(1,1)$-forms and can be thought of as metrics on $\Sigma$. But from the above equations we see that $a_1$ has strictly positive curvature while $a_2$ has strictly negative curvature. Therefore there are no compact solutions without singularities. The equation for an elliptic affine sphere in $\mathbb{R}^3$ \cite{lyz} appears as a special case of equation (\ref{d23toda}). Indeed the elliptic affine sphere equation is
\begin{equation}\label{tzitz1}
2(w_1)_{z\overline{z}} = -q\overline{q}e^{-4w_1} - e^{2w_1}.
\end{equation}
If we set $H_2 = q\overline{q}/H_1$ then equation (\ref{tzitz1}) yields a solution to (\ref{d23toda}) (strictly speaking it is a solution away from the zeros of $q$, however one can show the corresponding minimal immersion extends over the zeros of $q$). In fact we can explain this reduction more directly in terms of $G_2$ geometry. Equation (\ref{tzitz1}) corresponds to the equation for a cylindrical semi-flat $G_2$-manifold with scaling symmetry as in Section \ref{secredtosurf}.

There are other cases of minimal surfaces in $Q \subset \mathbb{R}^{3,3}$ that we can find similar equations for, namely real forms of the affine Toda equations for the affine Dynkin diagrams $A^{(1)}_1$ and $B^{(1)}_2$. These are minimal surfaces for which the image lies in a proper subspace of $\mathbb{R}^{3,3}$. We leave the details to the reader.


\addcontentsline{toc}{chapter}{Bibliography}

\end{document}